\documentclass[leqno, a4paper,11pt]{amsart}
\pdfoutput=1

\usepackage{enumitem}
\usepackage{amsmath,mathtools,mathrsfs,amsfonts,amssymb}
\usepackage[margin=2.4cm]{geometry}

 \mathchardef\ordinarycolon\mathcode`\:
\mathcode`\:=\string"8000
\begingroup \catcode`\:=\active
  \gdef:{\mathrel{\mathop\ordinarycolon}}
\endgroup

\usepackage{graphicx}
\usepackage{multirow}

\usepackage{hyperref}

\usepackage{setspace}
\makeatletter
\newcommand{\leqnomode}{\tagsleft@true\let\veqno\@@leqno}
\newcommand{\reqnomode}{\tagsleft@false\let\veqno\@@eqno}
\makeatother
\usepackage[english]{babel}

\newtheorem{thm}{Theorem}[section]

\newtheorem{lem}[thm]{Lemma}
\newtheorem{prop}[thm]{Proposition}
\theoremstyle{definition}
\newtheorem{defn}[thm]{Definition}
\newtheorem{rem}[thm]{Remark}
\newtheorem*{ex}{Example}
\numberwithin{equation}{section}

\newtheorem{fact}[thm]{Fact}

\newcommand{\bt}{\begin{thm}}
\newcommand{\et}{\end{thm}}
\newcommand{\bp}{\begin{prop}}
\newcommand{\ep}{\end{prop}}
\newcommand{\bd}{\begin{defn}}
\newcommand{\ed}{\end{defn}}
\newcommand{\bl}{\begin{lem}}
\newcommand{\el}{\end{lem}}
\newcommand{\bfa}{\begin{fact}}
\newcommand{\efa}{\end{fact}}
\newcommand{\bc}{\begin{corollary}}
\newcommand{\ec}{\end{corollary}}
\newcommand{\bex}{\begin{ex}}
\newcommand{\eex}{\end{ex}}
\newcommand{\ben}{\begin{enumerate}}
\newcommand{\een}{\end{enumerate}}

\newcommand{\ds}{\displaystyle}

\newcommand{\floor}[1]{\lfloor #1 \rfloor}

\newcommand{\sotto}[2]{#1_{#2}}
\newcommand{\lra}{\longrightarrow}
\newcommand{\rrr}{\rightarrow}
\newcommand{\ra}{\rightarrow}

\newcommand{\ideal}[1]{\sotto {{\mathcal I}}{#1}}

\newcommand{\exact}[3]
{0 \rrr #1 \rrr #2
\rrr #3 \rrr 0}

\newcommand{\pso}{\mathbb{P}^3}

\newcommand{\cac}{{\mathscr C}}

\newcommand{\coo}{{\mathcal O}}

\newcommand{\JJ}{{\mathcal J}}

\newcommand{\HH}{\mathrm{H}}

\newcommand{\rrrnk}{\mbox{rank}}

\makeatletter
\@namedef{subjclassname@2020}{
  \textup{2020} Mathematics Subject Classification}
\makeatother
\begin{document}

\title{Multiple Lines of Maximum Genus in $\pso$ }

\author{ Enrico Schlesinger}

\address{Dipartimento di Matematica, Politecnico di Milano \\ Piazza Leonardo da Vinci 32, 20133 Milano, Italy}

\email{enrico.schlesinger@polimi.it}
\subjclass{14H50, 14F06, 14F17}

\begin{abstract}
We introduce a notion of good cohomology for multiple lines in $\pso$ and we classify multiple lines
with  good cohomology up to multiplicity $4$. In particular, we show that the family of space curves of degree $d$,
not lying on a surface of degree $<d$, and of maximal arithmetic genus is not irreducible already for $d=4$ and $d=5$.
\end{abstract}
\thanks{The author was partially supported by   PRIN 2020 Squarefree Gr\"{o}bner  degenerations, special varieties and related topics, and is a member of GNSAGA of INDAM}

\keywords{Space curves, arithmetic genus, quasiprimitive multiple lines, Hilbert schemes}

\maketitle
\section{Introduction} \label{intro}
By a {\em space curve} we mean a locally Cohen-Macaulay purely one dimensional subscheme of $\pso$, the projective space over an algebraically closed field. Thus a curve is allowed to have several irreducible components and a nonreduced scheme structure, but it cannot have zero dimensional components - neither isolated nor embedded.
The most important invariants of a space curve $C$ are its {\em arithmetic genus} $g(C)=1-\chi \coo_C$, which does not depend on the embedding of $C$ in $\pso$;
its {\em degree} $\deg (C)$, which is defined through the Hilbert polynomial $\chi \left( \coo_C(n)\right)= n \deg (C) + 1 -g(C)$ and depends on the invertible sheaf $\coo_C(1)$, but not on the sections of  $\coo_C(1)$ that define the embedding of $C$ in $\pso$; and the {\em minimum degree} $s(C)$
of a surface that contains $C$, which {\em does depend} on the embedding in $\pso$.
The maximum genus problem for space curves asks to determine the most basic relation between these invariants, that is, what is the maximum
arithmetic genus $P(d,s)$ of a space curve of degree $d$ in $\pso$ that is not contained in a surface of degree $< s$ -  there is a huge literature on the maximum genus problem for {\em smooth}  space curves, but we will not be concerned with smooth curves in this paper.
The problem makes sense  for  pairs of integers $(d,s)$ satisfying $1 \leq s \leq d$  because there exists
a space curve of degree $d$ that is not contained in a surface of degree $< s$ if and only if $ d \geq s \geq 1$.

We now survey what is known about the maximum genus problem for (locally Cohen-Macaulay) space curves.
Beorchia \cite{beo2} proved a bound $B(d,s)$ for the maximum genus if the characteristic of the ground field is zero, and proved the bound is sharp if
$s \leq 4$; later the author \cite{sch-beo} gave a different proof of this bound valid in any characteristic.

\bt[\cite{beo2,sch-beo}]
\label{beo}
Let $C$ be a curve in $\pso$ of degree $d$ and genus $g$. Assume
that $C$ is not contained in any surface of degree $<s$.
Then $d \geq s$ and
\begin{equation}\label{Bbound}\tag{$\star$}
g \leq B(d,s) =
\begin{cases}
(s-1)d +1 - \binom{s+2}{3}, & \text{if~}  s \leq d \leq 2s, \\
\binom{d-s}{2}-
\binom{s-1}{3},& \text{if~} d \geq 2s+1.\\
\end{cases}
\end{equation}
\et
To prove sharpness one needs to construct, for each pair $d \geq s$, a curve of genus $B(d,s)$ not lying on a surface of degree $<s$.
The case $d=s$ is crucial because, if $P(s\!-\!1, s\!-\!1)=B(s\!-\!1,s\!-\!1)$,
it follows that $P(d,s)=B(d,s)$ for every
$d \geq 2s-1$ - see \cite{bls}. The aim of this paper is to propose a framework for the classification of curves achieving the maximum genus $B(d,d)$ in the basic case $\deg(C)=s(C)=d$, together with the computation of the first few relevant examples; even for $d=4$ this seems to be new.

One first observes that curves of degree $d$ are always contained in surfaces of degree $d$, and that those
that do not lie on a surface of degree $<d$ are supported
on either one line or two disjoint lines - see Proposition \ref{sequald} below.
It is therefore necessary to study curves supported on a line $L$ in $\pso$. We denote by $L_d$ the $(d\!-\!1)$-th neighbourhood  of $L$ in $\pso$: the ideal sheaf of $L_d$
is the $d$-th power  $\ideal{L}^d$ of the ideal sheaf of $L$. Any degree $d$ curve $C$ supported on $L$, for short a $d$-uple line, is contained in
$L_d$, and such a curve does not lie on a surface of degree $s<d$ if and only if  $\HH^0 \big(\ideal{C}(d\!-\!1)\big)= \HH^0 \big(\ideal{L}^d (d\!-\!1)\big)$
as the latter vector space vanishes. It is clear how to strengthen this requirement to deal with the fact that the support of a curve of maximum genus
in the case $d=s$ may consists of two disjoint lines, rather than only one.

\bd \label{cdl-def}
Fix integers $d \geq 1$ and $\ell \geq 0$.
We say that a degree $d$ curve $C$ supported on a line $L$ is a $C_{d,\ell}$ if
\begin{itemize}
  \item the genus of $C$ is
  $$
  g(C_{d,\ell})= B(d,d)-\ell \binom{d}{2} = -(d-1) - \binom{d}{3}-\ell \binom{d}{2}
  $$
  \item the only surfaces of degree $\ell +d-1$ containing $C$ are those containing the entire neighborhood $L_d$ of $L$ as well:
  $$\HH^0 \big(\ideal{C}(\ell +d-1)\big)= \HH^0 \big(\ideal{L}^d (\ell +d-1)\big).$$
\end{itemize}
\ed
\noindent
In particular, a $C_{d,0}$ is a $d$-uple line of genus $B(d,d)$ that does not lie on a surface of degree $<d$,
and the existence of a $C_{d,0}$ implies sharpness of Beorchia's bound $P(d,d)=B(d,d)$. But
the definition is tailored so that, for each  $1 \leq k \leq d-1$,  if $C$ and $D$ are respectively a $C_{k,d-k}$ and a $C_{d-k,k}$ whose supports are disjoint lines, then the  union of $C$ and $D$ is a also a curve of maximum genus, that is, a curve satisfying $\deg(C)=d$, $s(C)=d$ and $g(C)=B(d,d)$.

It was originally an idea of  Beorchia, see  \cite{bls}, that one should construct curves of maximum genus as $C_{d,0}$ for $d \equiv 2$ modulo $3$,
and adding a line to a $C_{d\!-\!1,1}$ when $d \equiv 0$ modulo $3$, or a suitable double line to a $C_{d\!-\!2,2}$ when $d \equiv 1$ modulo $3$.
Thus the problem of sharpness of the bound $B(d,d)$ is reduced to constructing $d$-lines with good cohomological properties when $d \equiv 2$ modulo $3$, and this  construction in \cite{bls} is reduced to an algebraic statement \cite[Conjecture B on p. 142]{bls}. Sammartano and the author \cite{ASES} are completing the proof of this Conjecture under the additional hypothesis the ground field has characteristic zero, thus showing the existence of curves $C_{d,\ell}$  for every $d \equiv 2$ modulo $3$ and proving sharpness of Beorchia's bound in the case $s=d$.

The main contribution of this paper is to show there are other components of curves of maximum genus $B(d,d)$ by giving examples of curves $C_{d,\ell}$ in cases $d=3$ and $d=4$. As a consequence,  we show that the fmaily of space curves satisfying  $\deg(C)=d$, $s(C)=d$ and $g(C)=B(d,d)$ is not irreducible and  contains curves that are scheme theoretically very different from the one constructed in  \cite{bls}. We hope this will be useful for the problem
of sharpness of Beorchia's bound in the intermediate range $s+1 \leq d \leq 2s$, as curves of maximum genus in that range have to be constructed adding a plane curve to a curve satisfying $\deg(C)=s$, $s(C)=s$ and $g(C)=B(s,s)$ \cite{sch-beo}.

Our main theorem is the classification of the curves $C_{d,\ell}$ for $d \leq 4$. For this we need the notion of {\em quasiprimitive multiple structure}
introduced in \cite{banica}, a notion that we review in section \ref{sectionqp}. A quasiprimitive $d$-line has an invariant, called type, that is a string of
$d\!-\!1$ integers $(a;b_2, \ldots, b_{d-1})$. A quasiprimitive $d$-line is primitive if $b_2= \ldots =b_{d-1}=0$, so that for a primitive $d$-line
the type is a single integer $a$.  It is trivial to note that a line is a $C_{1,\ell}$ for any $\ell \geq 0$, and a double line is a $C_{2,\ell}$ if and only if it has genus $-1-\ell$, or, equivalently, it is a primitive double line of type $a=\ell$. In section \ref{examples} we classify $C_{d,\ell}$'s for $d=3$ and $d=4$ proving
\newpage
\bt \label{main1} \

\begin{enumerate}
  \item A triple line is a $C_{3,\ell}$ if and only if it is quasiprimitive of type $(\ell;1)$.
The family of $C_{3,\ell}$ curves is irreducible of dimension $5 \ell +12$.
  \item A quadruple line is a $C_{4,\ell}$ if and only if it is a general quasiprimitive quadruple line of type $(\ell;2,2)$.
The family of $C_{4,\ell}$ curves is irreducible of dimension $9 \ell +21$.
\end{enumerate}
\et

Unfortunately, for $d \geq 5$ we don't have a classification. What we can say in general is that a $C_{d,\ell}$ is a quasiprimitive $d$-uple line
of type $(a;b_2, \ldots, b_{d-1})$ where $\ell \leq a \leq \ell + \floor{\frac{d-2}{3}}$. When $d \equiv 2$ modulo $3$, a strategy for constructing a primitive $C_{d,\ell}$ of type $a=\ell + \frac{d-2}{3}$ is proposed in \cite{bls}, and a proof that this works when the ground field has characteristic zero is being written up \cite{ASES}. But we don't know even for $d=5$  whether the quasiprimitive type is determined for a $C_{5,0}$ or whether the family of $C_{5,0}$'s is irreducible.

As an application of Theorem \ref{main1} we can show in the last section of the paper that the family of degree $d$ curves of maximum genus $B(d,d)$ that do not lie on a surface of degree $<d$ is not irreducible already for $d=4$ and $5$. Specifically

 \bt \label{main2} \

\begin{enumerate}
  \item The family of quadruple lines of maximum genus $B(4,4)=-7$ not lying on a cubic surface is not irreducible. It contains
  \begin{itemize}
    \item the $22$-dimensional irreducible family whose general member is the disjoint union of two double lines of genus $-3$;
    \item the $21$-dimensional family
      whose general member is the disjoint union of a line and a $C_{3,1}$;
  \end{itemize}
  the closures of these two families are different components of the Hilbert
      scheme $H_{4,-7}$ parametrizing space curve of degree $4$ and genus $-7$.
    \item The family of quintuple lines of maximum genus $B(5,5)=-14$ not lying on a quartic surface is not irreducible. It contains
     \begin{itemize}
       \item  the $30$-dimensional irreducible family whose general member is a general primitive quintuple line of type $a=1$;
       \item
          the $34$-dimensional family whose general member is the disjoint union of a line and a general $C_{4,1}$;
       \item
       the $35$-dimensional
          family whose general member is the disjoint union of a $C_{3,2}$ and a  $C_{2,3}$;
     \end{itemize}
     and there are no containment between the closures of these $3$ families in the Hilbert scheme $H_{5,-14}$.
\end{enumerate}
\et

\section*{Acknowledgements}
The Author would like to thank warmly Scott Nollet to whom he is indebted
for several discussions on the topic of multiple lines.

\section{Quasiprimitive multiple lines in $\pso$} \label{sectionqp}
By the term \begin{em} $d$-uple line\end{em} we will mean
a (locally Cohen-Macaulay) curve in $\pso$ that has degree $d$ and whose support is a line.
The notion of {\em quasiprimitive} multiplicity structure on a smooth curve was introduced by  Banica and Forster
\cite[$\S \; 3$]{banica}; we recall what it means in our context.

Let $C$ be a $d$-uple line with support $L$. Denote by $C_{j}$  the
subscheme of $C$ obtained by removing the embedded points from $C \cap L_j$ - as in the introduction, $L_j$ is the infinitesimal neighborhood of $L$ in $\pso$ defined by $\ideal{L}^j$ .  The {\it Cohen-Macaulay filtration} of
$C$ is:
\begin{equation}\label{filtration}
         L = C_{1} \subset C_{2} \subset \dots \subset C_{k} = C
\end{equation}
where $k$, $1 \leq k \leq d$, is the smallest integer such that $C \subset L_k$.
The quotients $\mathcal{L}_{j} = \ideal{C_{j}}/\ideal{C_{j+1}}$ are vector bundles on
$L$ and  $d= \deg (C) = 1 + \sum \rrrnk \; \mathcal{L}_{j}$.
The natural inclusions
$\ideal{C_{i}} \ideal{C_{j}} \subset \ideal{C_{i+j}}$ induce generically
surjective multiplication maps $\mathcal{L}_{i} \otimes \mathcal{L}_{j} \ra \mathcal{L}_{i+j}$ and
in particular we obtain generic surjections $\mathcal{L}_{1}^{j} \ra \mathcal{L}_{j}$.

A multiple line $C$ is {\it quasiprimitive} if it has generically embedding dimension two. This is the
case if and only if $\rrrnk \; \mathcal{L}_{1} = 1$, or, equivalently, $C$ does not contain the first infinitesimal
neighborhood $L_2$ of its support  $L$, so that the first filtrant $C_2$ has degree $2$ (and $C_j$ degree $j$ for each $j$).
If $C$ is quasiprimitive, then the generic
surjections $\mathcal{L}_{1}^{j} \ra \mathcal{L}_{j}$ of invertible sheaves yield effective divisors $D_{j}$
such that $\mathcal{L}_{j} \cong \mathcal{L}_{1}^{j}(D_{j})$; the multiplication maps
show that $D_{i} + D_{j} \leq D_{i+j}$.

For a quasiprimitive $d$-uple line $C$ in $\pso$, we define the {\em  type}
$\sigma(C)= (a;b_2,..,b_{d-1})$ of $C$ setting
$a= \deg (\mathcal{L}_1)$ and $b_j =\deg (D_j)$; it is convenient to set $b_1=0$ so that the inequalities
$$
b_i +b_j \leq b_{i+j}
$$
hold for every $i,j \geq 1$ such that $i+j \leq d-1$.

Finally, a $d$-uple line $C$ is called {\em primitive} if it has embedding dimension two everywhere. This is the case if and only if
$\mathcal{L}_{j} \cong \mathcal{L}_{1}^{j}$ for every $1 \leq j \leq d\!-\!1$, so $b_2= \ldots =b_{d-1}=0$ and the
type of $C$ simplifies to the single integer $a= \deg ( \mathcal{L}_{1})$.

\bp[Genus of a quasiprimitive multiple line]
Let $C$ be a quasiprimitive multiple line of type $(a;b_2,..,b_{d-1})$ in $\pso$. Then $a \geq -1$ and
\begin{equation}\label{genusCqpms}
g(C)= -(d-1) - \frac{a}{2}d(d-1) - \sum_{j=2}^{d-1} b_j 
\end{equation}
\ep

\begin{proof} Let $L$ be the support of $C$.
The inequality $a \geq -1$ follows from the fact that $\ideal{L,C_2} \cong \coo_L(a)$ is a quotient of the conormal
bundle $\ideal{L,L_2} \cong \coo_L(-1) \bigoplus \coo_L(-1)$. By definition of the type, $\ideal{C_{j},C_{j+1}}\cong \coo_L (ja+b_j)$. The formula for the genus
follows from the fact that $g(C)=\chi (\ideal{C})$.
\end{proof}

We next compute the dimension of the irreducible family of primitive $d$-uple lines of a given type $a \geq 0$.
Let $C$ be a primitive $d$-structure of type $a$ on the line $L$ in $\pso$. Given a subscheme $X \subset \pso$  we denote by the symbol
$\cac_{X}=\ideal{X}/\ideal{X}^2$ its conormal sheaf. Then \cite{banica} there exists an exact sequence
\begin{equation}\label{con1}
0 \lra \coo_L(da) \stackrel{\tau}{\lra}   \cac_{C} \otimes \coo_L \lra \cac_{L} \lra \coo_L\lra 0
\end{equation}
The morphism $\tau$ is induced by the inclusion $\ideal{L}^d \hookrightarrow \ideal{C}$ via the isomorphism
$\coo_L(de) \cong \ideal{L}^d/\ideal{L}^{d-1}\ideal{C_2}$. If $a \geq 0$, it follows that
\begin{equation}\label{conormalsplitting}
\cac_{C} \otimes \coo_L \cong  \coo_L(da) \oplus \coo_L(-a-2)
\end{equation}

By \cite[Proposition 2.3]{banica} the set of primitive $d+1$-structures $\tilde{C}$ that contain $C$ is parametrized by the
set of retractions $\beta:  \cac_{C} \otimes \coo_L \lra \coo_L(da)$ of $\tau$; the correspondence is given by
$\ideal{\tilde{C}}/\ideal{L}\ideal{C}= \mbox{Ker} (\beta)$. Therefore, if $a \geq 0$, the set of such $\tilde{C}$'s is parametrized
by the set of splittings of
$$
0 \lra \coo_L(da) \stackrel{\tau}{\lra}  \coo_L(da) \oplus \coo_L(-a-2) \lra \coo_L(-a-2) \lra 0
$$
hence by an affine space 
of dimension $(d+1)a+3$.
With a little extra effort one can check that the set $\mathcal{P}_L(d;a)$ of primitive $d$ structures on $L$ of type $a$ is an algebraic affine
bundle over  $\mathcal{P}_L(d\!-\!1;a)$, hence inductively that $\mathcal{P}_L(d;a)$ is a smooth variety of dimension
$$
(2a+3)+(3a+3)+ \cdots+(da+3)=\frac{a}{2}(d^2+d-2)+3(d-1)=(d-1)(3+\frac{a}{2}(d+2)).
$$
If we let the line $L$ vary as well, we obtain
\begin{equation}\label{dimensionPda}
\dim \mathcal{P}(d;a)= \frac{a}{2}(d-1)(d+2)+3d+1
\end{equation}
This is an interesting number as  primitive $d$-uple line are usually the generic point of a component of the Hilbert scheme parametrizing curves
of degree $d$ - see \cite{nthree} and \cite{ns-comp} for the first relevant examples.
When $d=2$ we recover the easy and well know fact that double lines of type $a \geq 0$, that is, of genus $-a-1 \leq -1$, form an irreducible family of dimension $2a+7$, which is a component of the Hilbert scheme if $a \geq 1$.

By a similar argument Nollet \cite[Corollary 2.6]{nthree} proves that the family $\mathcal{P}(3;a;b)$ of quasiprimitive triple lines of type $(a;b)$ with $a \geq 0$ is irreducible of dimension
\begin{equation}\label{dimensionP3ab}
\dim \mathcal{P}(3;a;b)= 5a+2b+10
\end{equation}
and one can prove, more generally, that the family $\mathcal{P}(d;a;0,\ldots,0,b)$ of quasiprimitive $d$-lines of type $(a;0,\ldots,0,b)$,
with $a \geq 0$,
 is irreducible of dimension
\begin{equation}\label{dimensionPda0b}
\dim \mathcal{P}(d;a;0,\ldots,0,b)= \frac{a}{2}(d-1)(d+2) +3d+2b+1.
\end{equation}

With some extra effort  Nollet and the author  \cite[Proposition 2.3]{ns-comp} prove that the family  $\mathcal{P}(d;a;b,c)$ of quasiprimitive
 $4$-lines of type $(a;b,c)$, with $a \geq 0$, is irreducible of dimension
  \begin{equation}\label{dimensionP4abc}
\dim \mathcal{P}(d;a;b,c)= 9a+2b+2c+13
\end{equation}
The extra effort goes into proving \cite[Lemma 2.2]{ns-comp} that, for a quasiprimitive triple line $C$ of type $(a;b)$ with $a\geq 0$ supported on the line $L$,
the restriction of the conormal sheaf $\cac_{C} \otimes \coo_L$ has torsion, and modulo torsion is isomorphic to
$ \coo_L(3a+b) \oplus \coo_L(-a-b-2)$. A similar argument works for a quasiprimitive $d$-line $C$ of type $(a;0,\ldots,0,b)$,
with $a \geq 0$ and shows $\cac_{C} \otimes \coo_L$  modulo torsion is isomorphic to
$ \coo_L(da+b) \oplus \coo_L(-a-b-2)$. We can then prove the irreducibility of the family of quasiprimitive $d$-uple lines
of type $(a;0,\ldots,0,b,c)$, with $a \geq 0$, and compute its dimension:
 \begin{equation}\label{dimensionPda0bc}
\dim \mathcal{P}(d;a;0,\ldots,0,b,c)= \dim \mathcal{P}(d\!-\!1;a;0,\ldots,0,b)+ da+2c+3
\end{equation}

\section{Rough classification of curves with $s(C)= \deg (C)$}

The following proposition is easy and certainly well-known, see \cite[Remark 6.8]{sch-tran}, but we include a proof for the convenience of the reader:
\bp \label{sequald}
For a space curve $C$, the inequality  $s(C) \leq \deg (C)$  holds; if $s(C)=\deg (C)$, then
\begin{enumerate}
\item every subcurve $D \subseteq C$ also satisfies $s(D)=\deg (D)$;

\item the curve $C_{red}$ is either a line or the disjoint union of two lines, and on each
line in its support $C$ has the structure of a quasiprimitive multiple line satisfying $a \geq 0$
(here $a=\deg (\mathcal{L}_{1})$ is the first integer appearing in the type of $C$).
\end{enumerate}

\ep

\begin{proof}
The inequality $s(C) \leq \deg (C)$ is proven for example in \cite{sch-beo}. Suppose from now on that
$s(C)=\deg(C)$. If $D \subseteq C$ and $S$ is a surface of degree $s(D)$ containing $D$, there is an exact sequence
$$\exact{\ideal{E}(-s(D))}{\ideal{C}}{\ideal{C \cap S, S}}$$
where $E$, the subscheme of $C$ residual to $C \cap S$, is a locally Cohen-Macaulay curve of degree
$$\deg(E)=\deg (C)-\deg (C \cap S) \leq \deg (C)- \deg(D).$$
Hence $s(E) \leq  \deg (C)- \deg(D)$. On the other hand, by the above exact sequence $ s(C) \leq s(E)+s(D)$, therefore
$ s(C) \leq  \deg (C)- \deg(D)+s(D)$.
This together, with the hypothesis $s(C)=\deg(C)$, implies $\deg(D) \leq s(D)$, hence equality must hold.

For an irreducible and reduced curve $D$, the equality $s(D)= \deg(D)$ can hold only if $D$ is a line - cf. \cite[Proposition 3.2]{sch-beo}. Note that $C$ cannot contain
the union $D$ of two lines meeting at one point because such a $D$ has $s(D)=1<\deg(D)$.
Thus, the support of $C$ is a union of disjoint lines. Since the union of $3$ disjoint lines
lies on a quadric surface, the support of $C$ consists of at most two lines. Since the first infinitesimal neighborhood of a line in $\pso$
has degree $3$ and is contained in a quadric surface, it cannot be contained in $C$, so $C$ has a
quasiprimitive structure on each line in its support. Finally, since any degree $2$ subcurve of $C$ is not contained in a plane while a double line of type $-1$ is contained in a plane, we must have $a \geq 0$.
\end{proof}

%
%

\section{Multiple lines with good cohomological properties}

Fix homogeneous coordinates $x,y,z,w$ on $\pso$ so that $L$ is the line of equations
$x=y=0$. The projection $\pi: L_d \ra L$ from the line $M$ of equations $z=w=0$ corresponds
to the inclusion of coordinate rings
$$H^{0}_* (\coo_{L}) = k[z,w] \cong k[x,y,z,w]/(x,y) \hookrightarrow k[x,y,z,w]/(x,y)^d \cong H^{0}_* (\coo_{L_d})$$
From the isomorphism of $k[z,w]$-modules
\begin{equation}\label{RingLd}
H^{0}_* (\coo_{L_d})= \frac{k[x,y,z,w]}{(x^d, x^{d-1}y, \ldots, y^d)} \cong \bigoplus_{i=0}^{d-1} \left(k[z,w](-i)\right)^{\oplus(i+1)}
\end{equation}
it follows
\begin{equation}\label{spectrum Ld}
\pi_* \coo_{L_d} \cong
\bigoplus_{i=0}^{d-1} \left(\coo_L(-i)\right)^{\oplus(i+1)}.
\end{equation}

\bp \label{Cdabyvanish}
Let $C$ be a curve of degree $d$ supported on the line $L$, and let $\pi: L_d \rightarrow L$ denote the projection from
a line $M$ disjoint from $L$. Fix an integer $\ell \geq 0$. The following conditions are equivalent
\begin{enumerate}
\item the genus of $C$ is
$$g(C)=-(d-1) - \binom{d}{3}-\ell \binom{d}{2}=B(d,d)-\ell \binom{d}{2}$$
and
$$\HH^0 \big(\ideal{C}(\ell +d-1)\big)= \HH^0 \big(\ideal{L}^d (\ell +d-1)\big).$$
\item the genus of $C$ is
$g(C)=B(d,d)-\ell \binom{d}{2}$
and $\HH^1 \big(\ideal{C}(\ell +d-1)\big)= 0$.
\item
the sheaf $\pi_* \ideal{C,L_d}$ is isomorphic to $
\left(
\coo_L (-d-\ell)
\right)
^{\oplus \frac{d(d-1)}{2}}$.
\end{enumerate}
\ep

\begin{proof}
If $C$ has degree $d$ and genus $B(d,d)-\ell \binom{d}{2}=g(C_{d,\ell})$, then the rank and the degree of the locally free
sheaf $\mathcal{E}=\pi_* \ideal{C,L_d}$ are the same
as those of $\coo_L (-d-\ell) ^{\oplus \frac{d(d-1)}{2}}$. Thus  $\mathcal{E} \cong \coo_L (-d-\ell) ^{\oplus \frac{d(d-1)}{2}}$
if and only if $h^0 \mathcal{E} (\ell+d-1)=0$, which is equivalent to $h^1 \mathcal{E} (\ell+d-1)=0$ because $\chi \mathcal{E} (\ell+d-1)=0$.
As $h^1 \ideal{L_d} (\ell+d-1)  = h^2 \ideal{L_d} (\ell+d-1)=0$, these vanishings are equivalent to those in the statement.
\end{proof}

\bd
Given a pair of integers $d \geq 1$ and $\ell \geq 0$, we say that a $d$-line is a $C_{d,\ell}$ if it satisfies the equivalent conditions of proposition \ref{Cdabyvanish}.
\ed
\noindent
Note that a $C_{d,0}$ is a curve of degree $d$, not lying on a surface of degree $<d$,
of maximum genus $B(d,d)$: if a $C_{d,0}$ exists for a given $d$, then Beorchia's bound $B(d,d)$ is sharp.

A line is a $C_{1,\ell}$ for every $\ell$. A  double line of genus $-\ell-1 <0$ is a $C_{2,\ell}$, and conversely.
Indeed, all non planar double lines arise as follows - see for example \cite{nthree}: take a smooth
surface $S$ containing $L$ of degree $\ell+2 \geq 2$ and let $C$ be the divisor $2L$ on $S$.
By adjunction
$$\ideal{L,C} \cong \coo_S (-L) \otimes \coo_L \cong \coo_L (\deg (S)-2))= \coo_L(\ell).$$
From the exact sequence
$$
\exact{\ideal{C,L_2}}{\ideal{L,L_2}}{\ideal{L,C} \cong \coo_L(\ell)}
$$
we see
$$\ideal{C,L_2} \cong \coo_L (-\ell-2),$$
that is,
$C$ is a $C_{2,\ell}$. Since $\ideal{L,C} \cong \coo_L(\ell)$, the double
line $C$ is primitive of type $a=\ell \geq 0$.

To analyze higher degree cases, we introduce an intermediate notion.
\bd
Given integers $d \geq 1$ and $\ell \geq 0$, we say that a $d$-uple line $C$ satisfies condition $\mathscr{C}_{d,l}$ if
\begin{equation}\leqnomode\label{conditiondl}
 \hspace{2cm} h^0 (\ideal{C,L_d}(\ell + d -1))=0 \tag{$\mathscr{C}_{d,l}$}
\end{equation}
or, equivalently, $
\HH^0 \big(\ideal{C}(n)\big)= \HH^0 \big(\ideal{L}^d (n)\big)$ for every $n \leq \ell+d-1$.
\ed
For $\ell=0$, the condition $\mathscr{C}_{d,0}$ means simply that $C$ is not contained in any surface of degree $<d$, but for $\ell>0$ the condition $\mathscr{C}_{d,\ell}$ is stronger as it says
that the only surface containing $C$ of degree up to $\ell + d -1$ are those containing the whole infinitesimal neighborhood $L_d$.
By Proposition  \ref{maxgenusld} below a $C_{d,\ell}$ is a $d$-uple line of maximum genus among those satisfying condition $\mathscr{C}_{d,l}$.
Note that a double line $C$ of type $a$, that is, of genus $-a-1$, satisfies condition $\mathscr{C}_{2,l}$ if and only if $a \geq \ell$
because $\ideal{C,L_2} \cong \coo_L (-a-2)$.

\bl \label{cdl}
If $C$ satisfies condition $\mathscr{C}_{d,l}$ and $ D \subset C$ is a locally Cohen-Macaulay subcurve of degree $k$,
then $D$ satisfies condition $\mathscr{C}_{k,l}$.
\el

\begin{proof}
This follows from $I_L^{d-k} \ideal{D} \subseteq \ideal{C}$.
\end{proof}

\begin{rem}
Unfortunately, it is not true that a degree $k$ subcurve $D$ of a $C_{d,\ell}$ is a $C_{k,\ell}$;
the point of the previous lemma is that at least $D$ satisfies condition $\mathscr{C}_{k,l}$.
\end{rem}

A $C_{d,\ell}$, assuming it exists, has maximum genus among degree $d$ multiple lines whose ideal agrees with that of $L_d$ up to degree $\ell+d-1$:
\bp \label{maxgenusld}
Suppose $C$ is $d$-uple line with support $L$.
If $\ell \geq 0$ and
\[
\HH^0 \big(\ideal{C}(n)\big)= \HH^0 \big(\ideal{L}^d (n)\big)\quad \text{for~} n \leq \ell+d-1,
\]
then $\ds g(C) \leq  -(d-1) - \binom{d}{3}-\ell \binom{d}{2}$.
\ep

\begin{proof}
By hypothesis, $\HH^0 \big(\ideal{C}(d\!-\!1)\big)=0$, hence $C$ is a curve of degree $d$ that does not lie on a surface of degree $<d$.
It follows $\HH^1 \big(\coo_{C}(n)\big)=0$ for $n \geq -1$ by \cite[Proposition 3.2]{sch-beo}. Hence
$$
d(\ell+d-1)+1-g(C)= h^0 (\coo_C(\ell+d-1)) \geq h^0 (  \coo_{\pso} (\ell+d-1))-h^0 (\ideal{L}^d (\ell+d-1))
$$
which is equivalent to $g(C) \leq g(C_{d,\ell})$
because
$$
d(\ell+d-1)+1-g(C_{d,\ell})= h^0 (\coo_{\pso} (\ell+d-1))-h^0 (\ideal{L}^d (\ell+d-1)).
$$
\end{proof}
In \cite[p. 141]{bls}, with a different terminology, it is noted that one could prove sharpness of Beorchia's bound in the case $d=s$ by constructing curves $C_{d,\ell}$ for all
 $d \equiv 2$ modulo $3$ and $\ell=0,1,2$: indeed  for $d \equiv 2$ modulo $3$, a $C_{d,0}$ has genus $B(d,d)$; when $d \equiv 0$ modulo $3$  the disjoint union of  a line and a $C_{d\!-\!1,1}$ has genus $B(d,d)$; finally, when $d \equiv 1$ modulo $3$ the disjoint union of a $C_{d\!-\!2,2}$
  and a double line of genus $1\!-\!d$ has genus $B(d,d)$. We introduced the notion of a $C_{d,\ell}$ to formalize and generalize this remark as follows:

\bp \label{CunionC}
Suppose $1 \leq k \leq d-1$ and $C$ and $D$ are respectively a $C_{k,d-k}$ and $C_{d-k,k}$ whose supports are disjoint.
Then the disjoint union of $C$ and $D$ is a curve of degree $d$ and genus $B(d,d)$ that does not lie on a surface
of degree $d-1$.
\ep
\begin{proof}
A direct calculation shows
$$
g (C \cup D)= g (C_{k,d-k}) + g (C_{d-k,k})-1 = B(d,d).
$$
Thus we only need to show that $C \cup D$ is not contained in a surface of degree $<d$. By way of contradiction,
suppose $F$ is the equation of a degree $d\!-\!1$ surface $S$ containing $C\cup D$. We can assume the support
of $C$ is the line of equations $x=y=0$ and the support of $D$ is the line $z=w=0$. By assumption, the polynomial $F$
must lie in $(x,y)^{k}$ because $C \subset S$ and in $(z,w)^{d-k}$  because $D \subset S$, but this contradicts $\deg (F)=d-1$.
\end{proof}

A $d$-uple line $C$ that is a $C_{d,\ell}$ is quasiprimitive, and
there are some obvious numerical constraints on the type of $C$:

\bp \label{num-constr}
Suppose $d \geq 2$ and $C$ is a $C_{d,\ell}$. Then $C$ is quasiprimitive.
If the type of $C$ is $(a;b_2,\ldots,b_{d-1})$, then $\ell \leq a \leq \ell + \floor{\frac{d-2}{3}}$
and
$$\sum_{j=2}^{d-1} b_j + (a-\ell) \binom{d}{2}= \binom{d}{3}
$$
In particular, if $C$ is primitive, then $d \equiv 2$ modulo $3$ and $a=  \ell + \frac{d-2}{3}$.
\ep
\begin{proof}
If $d$-uple line $C$ is a $C_{d,\ell}$, then in particular it does not lie on a surface of degree $< d$, hence it is quasiprimitive.
By Lemma \ref{cdl} the double line $C_2$ contained in $C$ satisfies condition $\mathscr{C}_{2,\ell}$ and has type $a$, hence
$a \geq \ell$.

Comparing the genus of a $C_{d,\ell}$ with the formula for the genus of a quasiprimitive multiple line, we obtain the equality:
$$
a  \binom{d}{2} + \sum_{j=2}^{d-1} b_j = \ell \binom{d}{2} + \binom{d}{3}
$$
If $C$ is primitive, that is, all the $b_j$'s are zero, it follows $d \equiv 2$ modulo $3$ and $a=  \ell + \frac{d-2}{3}$.
For an arbitrary  $C$, the integers $b_j$'s are nonnegative, and the above equality implies $a \leq \ell + \floor{\frac{d-2}{3}}$.
  \end{proof}

%

\section{Examples of low degree} \label{examples}

Triple lines have been studied by Nollet \cite{nthree}. In particular, he shows that the set of quasiprimitive triple
lines of type $(a;b)$ with $a,b \geq 0$ is nonempty and irreducible of dimension $5a+2b+10$ \cite[Corollary 2.6]{nthree}.

\bp[Triple lines] \label{3lines}
Fix an integer $\ell \geq 0$.
A triple line $C$ satisfies condition $\mathscr{C}_{3,l}$ if and only if it is quasiprimitive of type $(a;b)$
and either $a = \ell$ and $b \geq 1$, or $a \geq \ell + 1$.
Furthermore, $C$ is a $C_{3,\ell}$ if and only if $a=\ell$ and $b=1$. In particular,
the family of $C_{3,\ell}$ curves is irreducible of dimension $5 \ell +12$.
\ep

\begin{proof}
If $C$ satisfies condition $\mathscr{C}_{3,l}$, then $C$ is quasiprimitive because it does not lie on a surface of degree $2$.
Suppose the type of $C$ is $(a;b)$.
Then $C$ contains a unique double line $C_2$,
$
\ideal{L,C_2} \cong \coo_L (a)$, $\ideal{C_2,C} \cong \coo_L(2a+b)$.
By Lemma \ref{cdl} the double line $C_2$ satisfies condition $\mathscr{C}_{2,l}$, therefore $a \geq \ell$.
Consider the exact sequence of $\coo_L$-modules
$$
\exact{\frac{\ideal{C}}{\ideal{L}\ideal{C_2}}}{\frac{\ideal{C_2}}{\ideal{L}\ideal{C_2}}}{\coo_L(2a+b)}.
$$
As $\frac{\ideal{C_2}}{\ideal{L}\ideal{C_2}} \cong \coo_L(2a) \oplus \coo_L(-2-a)$, we conclude
$
\frac{\ideal{C}}{\ideal{L}\ideal{C_2}}\cong \coo_L(-2-a-b)
$. Then note that there is an obviously surjective
map of $\coo_L$-modules
$$
\alpha: \ideal{L,L_2} \otimes \ideal{C_2,L_2} \ra \frac{\ideal{L}\ideal{C_2}}{\ideal{L}^3}
$$
The sheaf on the left is isomorphic to $\coo_L (-3-a)^{\oplus 2}$, while the sheaf on the right is locally free of rank two,
thus $\alpha$ must be an isomorphism and $\frac{\ideal{L}\ideal{C_2}}{\ideal{L}^3} \cong \coo_L (-3-a)^{\oplus 2}$.
Finally, from the exact sequence
$$
\exact{\frac{\ideal{L}\ideal{C_2}}{\ideal{L}^3}}{\ideal{C,L_3}}{\frac{\ideal{C}}{\ideal{L}\ideal{C_2}}}
$$
we conclude $\pi_* \ideal{C,L_3} \cong \coo_L(-3-a)^{\oplus 2} \oplus \coo_L(-2-a-b)$ so that $C$
satisfies condition $\mathscr{C}_{3,l}$ if and only if either $a=\ell$ and $b \geq 1$ or
$a \geq \ell +1$, and $C$ is a  $C_{3,\ell}$ if and only if it has type  $(\ell;1)$.
\end{proof}

The case of quadruple lines is more difficult because the type of a quasiprimitive quadruple line $C$ does not determine its postulation,
that is, the sequence $n \mapsto h^0 (\ideal{C}(n))$. We show that a quadruple line is a  $C_{4,\ell}$ if and only if it is a
{\em sufficiently general} quasiprimitive quadruple line of  type $(\ell;2,2)$.
\bt[Quadruple lines] \label{4lines}
If a quadruple line $C$ satisfies condition $\mathscr{C}_{4,l}$, then $C$ is quasiprimitive of type $(a;b_2,b_3)$ and either
$a= \ell$ and $b_2 \geq 2$, or $a \geq \ell +1$.

Furthermore, if a quadruple line $C$ is a $C_{4,\ell}$, then $C$ is quasiprimitive of type $(\ell;2,2)$; and a
sufficiently general quasiprimitive quadruple line of  type $(\ell;2,2)$ is a $C_{4,\ell}$. Such curves form a nonempty
irreducible family of dimension $9 \ell +21$.
\et

\begin{proof}
Quadruple lines have been studied in \cite{ns-comp}. In particular, the dimension of the family of quasiprimitive quadruple
lines of type $(a;b,c)$ is computed in \cite[Proposition 2.3]{ns-comp}) and we summarized the argument on page \pageref{dimensionP4abc}.
Let $C$ be a quasiprimitive $4$-uple line
of type $(a;b_2,b_3)$. Then $C$ contains a unique double line $C_2$ and a unique triple line $C_3$, and
$$
\ideal{L,C_2} \cong \coo_L (a), \quad \ideal{C_2,C_3} \cong \coo_L(2a+b_2), \quad \ideal{C_3,C} \cong \coo_L(3a+b_3).
$$
If $C$ satisfies condition $\mathscr{C}_{4,l}$, then $C$ is quasiprimitive
and $C_3$ satisfies condition $\mathscr{C}_{3,l}$, hence by \ref{3lines} either $a=\ell$ and $b_2 \geq 1$ or $a \geq \ell +1$.

It remains to exclude the case $a=\ell$ and $b_2=1$. By \cite[Lemma 2.2]{ns-comp}, if one defines $\JJ= \ideal{L}\ideal{C_3}+\ideal{C_2}^2$, then
$\ideal{C_3}/\JJ \cong \coo_L (3a+b_2) \oplus \coo_L (-a-b_2-2)$. When $b_2=1$, it follows $h^0 ((\ideal{C_3}/\JJ) (a+3))= 4a+6$.
On the other hand, by the proof of \ref{3lines},
$
\pi_* \ideal{C_3,L_3} \cong \coo_L(-3-a)^{\oplus 3}
$
hence
$$
h^0 (\ideal{C_3} (a+3)) =  h^0 (\ideal{L_4} (a+3)) + 4(a+1)+3.
$$
Now suppose $a=\ell$ and let $p=d+\ell-1=a+3$: then
$$
h^0 (\ideal{C}(p)) \geq h^0 (\JJ (p)) \geq h^0 (\ideal{C_3} (p))- h^0 ((\ideal{C_3}/\JJ) (p))=
h^0 (\ideal{L_4} (p))+1.
$$
We conclude that $C$ does not satisfy condition $\mathscr{C}_{4,l}$ when $a=\ell$ and $b_2=1$.

Now suppose $C$ is a $C_{4,\ell}$, then looking at the genus $C$ we see $C$ must
be  quasiprimitive of type $(\ell;2,2)$. Let us show that a sufficiently general quasiprimitive quadruple line $C$ of  type $(\ell;2,2)$ is a $C_{4,\ell}$.

Using the same notation as above, in this case $\ideal{C_3,C} \cong \coo_L(3a+2)$ and
$\ideal{C_3}/\JJ \cong \coo_L (3a+2) \oplus \coo_L (-a-4)$, hence $\ideal{C}/\JJ \cong \coo_L (-a-4)$,
thus it is enough to show  that, for a general choice of
$C_3$, the sheaf
$\mathcal{F}=
\pi_* \JJ/\ideal{L}^4$ is isomorphic to $ \coo_L (-a-4)^{\oplus 5}$.
To prove this, we recall Nollet's description \cite{nthree} of the ideal of a quasiprimitive $3$-line $C_3$ of type $(a;2)$:
there exist forms
\begin{itemize}
  \item $f,g\in k[z,w]$ of degree $a+1$ with no common zero;
  \item $p,q \in k[z,w]$ of degree $2$ and $3a+4$ respectively, with no common zero;
  \item $r,s,t \in k[z,w]$ of degree $a+2$ such that $q=rf^2+sfg+tg^2$;
\end{itemize}
such that
\begin{enumerate}
  \item[i)] in the exact sequence
  $$
  \exact{\frac{\ideal{C_3}}{\ideal{L}\ideal{C_2}}}{\frac{\ideal{C_2}}{\ideal{L}\ideal{C_2}} \cong
  \coo_L (2a) \oplus \coo_L (-a\!-\!2)}{\ideal{C_2,C_3} \cong \coo_L (2a+2)}
$$
the last map is given by $[p, q]$;
  \item[ii)]   the homogeneous ideal of $C_3$ is
  $$
  I_{C_3}=I_L^3 + <xF,yF,G>
  $$
  where $F=xg-yf$ and $G=pF-rx^2-sxy-ty^2$.
\end{enumerate}
Let  $H_1 = xG$ and $H_2= yG$, and
consider the map $\beta$ of $\coo_L$-modules
$$
\mathcal{E}=\coo_L(-2a\!-\!4) \oplus \coo_L (-a\!-\!5)^{\oplus 2} \oplus \coo_L (-a\!-\!4)^{\oplus 3}
\stackrel{\beta}{\longrightarrow}
 \mathcal{G}=\pi_* \frac{\ideal{C_3}}{\ideal{L}^4}
$$
that sends the generators of $\mathcal{E}$ to the classes of $F^2$, $H_1$, $H_2$, $x^2F$,$xyF$, $y^2F$; as all these polynomials
are in $I_C^2+I_L I_{C_3}$, the map $\beta$ factors through the inclusion $\mathcal{F}=
\pi_* \JJ/\ideal{L}^4\hookrightarrow \mathcal{G}$.
By the proof of Proposition \ref{3lines} $\pi_* \frac{\ideal{C_3}}{\ideal{L}^3}$ is isomorphic
to $\coo_L(-a-3)^{\oplus 2} \oplus \coo_L(-a-4)$, and so
$$
\mathcal{G}=\pi_* \frac{\ideal{C_3}}{\ideal{L}^4} \cong
\coo_L(-3)^{\oplus 4} \oplus \coo_L(-a\!-\!3)^{\oplus 2} \oplus \coo_L(-a\!-\!4)
$$
with generators corresponding to $x^3$, $x^2y$, $xy^2$, $y^3$, $xF$, $yF$ and $G$. With respect to the chosen basis, the matrix of $\beta: \mathcal{E} \ra \mathcal{G}$ is
$$
\begin{bmatrix}
  0 & -r & 0 & g & 0 & 0 \\
  0 & -s & -r & -f & g & 0 \\
  0 & -t & -s & 0 & -f & g \\
  0 & 0 & -t & 0 & 0 & -f\\
  g & p & 0 & 0 & 0 & 0 \\
  -f & 0 & p & 0 & 0 & 0 \\
  0 & 0 & 0 & 0 & 0 & 0
\end{bmatrix}
$$
Generically, the map $\beta$ has rank $5$ as one can see for example by computing the $5 \times  5 $ minors of its matrix. On the other hand,
the section
$[p,-g,f,-r,-s,-t]^T$ of  $H^0({\mathcal{E}(2a+6)})$ is in the kernel of $H^0 (\beta(2a+6))$ by a direct check or  because
$$
pF^2-gH_1+fH_2-r x^2 F-sxyF-ty^2F=0.
$$
As $f$ and $g$ have no common zeros on $L$, we conclude that the kernel of $\beta$ is isomorphic to $\coo_L(-2a-6)$. As we have already observed, the image of $\beta$ is contained in $\mathcal{F}= \pi_* \JJ/\ideal{L}^4$
hence we have an exact sequence
$$
0 \ra \coo_{L}(-2a-6) \ra \mathcal{E} \stackrel{\beta}{\ra}\mathcal{F}
$$
Finally, from the exact sequence $\displaystyle \exact{\frac{\JJ}{\ideal{L}^4}}{\frac{\ideal{C_{3}}}{\ideal{L}^4}}{\frac{\ideal{C_3}}{\JJ}}$
we compute that $\mathcal{F}=\pi_* \JJ/\ideal{L}^4$ has the same rank $5$ and the same degree $-5a-20$ as the image of $\beta$,
hence  $ \mathcal{E} \stackrel{\beta}{\ra}\mathcal{F}$ is surjective.

%
We can now show  that $\mathcal{F} \cong \coo_L (-a-4)^{\oplus 5}$ if $C_3$ is general. As $\mathcal{F}$ is locally free of the same rank and degree
as $\coo_L (-a-4)^{\oplus 5}$, it is enough to prove $H^0 \mathcal{F}^\vee (-a-5) =0$.
Dualizing and twisting the above exact sequence we obtain
$$
\exact{\mathcal{F}^\vee (-a-5)}{\coo_L(-1)^{\oplus 3} \oplus \coo_L^{\oplus 2} \oplus \coo_L (a-1)}{\coo_L(a+1)}
$$
Thus what we need is that the map
$$
H^0 \coo_L^{\oplus 2} \oplus H^0 \coo_L (a-1) \stackrel{[f,-g,p]}{\ra} H^0 \coo_L(a+1)
$$
be injective. Now this is certainly the case is $f$, $g$ and $p$ are chosen general, as it is injective if we choose
$f=z^{a+1}$, $g=w^{a+1}$ and $p=zw$.
\end{proof}

\section{Families of maximum genus curves of low degree}
To summarize, the curves $C_{d,\ell}$ of which we know existence are:
\begin{enumerate}
\item when $d=3$, quasiprimitive multiple lines of type $(\ell;1)$ - this paper, Proposition \ref{3lines};
\item when $d=4$, quasiprimitive multiple lines of type $(\ell;2,2)$ - this paper, Theorem \ref{4lines};
\item when $d=3m-1 \leq 119$, primitive multiple lines of type $a= \ell +\frac{d-2}{3}$ - these are constructed in \cite{bls}, with the aid of {\em Macaulay2} for $m \geq 4$; at least in characteristic zero, in \cite{ASES} we will show how to extend this result to all degrees $d \equiv 2$ modulo $3$.
\end{enumerate}
One would be tempted to guess from these examples that the family  of $C_{d,\ell}$'s supported on a line $L$ is irreducible, but there might be counterexamples already for $d=5$:  I don't know if there are quasiprimitive quintuple lines of type $(0;2,2,6)$ or of type $(0;2,3,5)$ that do not lie on a quartic surface, or whether, if they exist, they lie in the closure of the family of primitive $5$ lines of type $a=1$.

We close the paper enumerating the known examples of degree $d$ curves of maximum genus $B(d,d)$ not lying on surfaces of degree $<d$ for
small $d$, thereby proving Theorem \ref{main2} in the introduction.

For $d=2$,  the maximum genus $B(2,2)$ is $-1$ and every curve of degree $2$ and $g=-1$ is not contained in a plane;  the $7$ dimensional irreducible family of double lines of genus $-1$,
that is of $C_{2,0}$'s, is contained in the closure of the family of two disjoint lines, which is the general member of the $8$ dimensional Hilbert scheme $H_{2,-1}$.

For $d=3$,  the maximum genus $B(3,3)$ is $-3$; Proposition \ref{CunionC} provides two irreducible families of degree $3$ curves of maximum genus
$-3$ not lying on a quadric: the $12$ dimensional family of quasiprimitive triple lines of type $(0;1)$, and the $13$ dimensional family whose general member is the disjoint union of a line and a double line of genus $-2$; the first family is in the closure of the second by \cite[Proposition 3.3]{nthree}.

For $d=4$, the family of quadruple lines of maximum genus $B(4,4)=-7$ not lying on a cubic surface is not irreducible.
It contains by Proposition \ref{CunionC}
  \begin{itemize}
    \item the $22$-dimensional irreducible family $F_1$ whose general member is the disjoint union of two double lines of genus $-3$;
    \item the $21$-dimensional family $F_2$ whose general member is the disjoint union of a line and a quasiprimitive triple line of type $(1;1)$;
     \item the $21$-dimensional family $F_3$ whose general member is a general quasiprimitive quadruple line of type $(0;2,2)$.
  \end{itemize}
It is clear the second family $F_2$ cannot be in the closure of family $F_1$, and the closure of these two families are in fact a component of the Hilbert scheme $H_{4,-7}$ by \cite[Theorem 6.2]{ns-comp}, while the family $F_3$ is in the closure of $F_1$ by  \cite[Proposition 3.3]{ns-comp}.

For $d=5$,  the family of quintuple lines of maximum genus $B(5,5)=-14$ not lying on a quartic surface is not irreducible. It contains
     \begin{itemize}
           \item
       the $35$-dimensional irreducible
          family $G_1$ whose general member is the disjoint union of a $C_{3,2}$ and a  $C_{2,3}$;
          \item
          the $34$-dimensional irreducible family $G_2$ whose general member is the disjoint union of a line and a general $C_{4,1}$;
       \item  the $30$-dimensional irreducible family $G_3$ whose general member is a general primitive quintuple line of type $a=1$.
     \end{itemize}
     and there are no containment between the closures of these $3$ families in the Hilbert scheme $H_{5,-14}$: one reason for which $G_3$ is not in the closure of $G_1$ or $G_2$ is that its general member is a curve that has embedding dimension two at each of its points, a property that does not hold for any curve in $G_1$ and $G_2$.

\end{document}